\documentclass[leqno,12pt]{article}

\usepackage{amssymb}
\usepackage{euscript}
\usepackage[dvips]{graphicx}
\usepackage{flafter}
\usepackage{pstricks}

\usepackage{mathrsfs} 

\evensidemargin\oddsidemargin

\begin{document}

\baselineskip=18pt
\setcounter{page}{1}

\renewcommand{\theequation}{\thesection.\arabic{equation}}
\newtheorem{theorem}{Theorem}[section]
\newtheorem{lemma}[theorem]{Lemma}
\newtheorem{proposition}[theorem]{Proposition}
\newtheorem{corollary}[theorem]{Corollary}
\newtheorem{remark}[theorem]{Remark}
\newtheorem{fact}[theorem]{Fact}
\newtheorem{problem}[theorem]{Problem}
\newtheorem{example}[theorem]{Example}

\newcommand{\eqnsection}{
\renewcommand{\theequation}{\thesection.\arabic{equation}}
    \makeatletter
    \csname  @addtoreset\endcsname{equation}{section}
    \makeatother}
\eqnsection

\def\r{{\mathbb R}}
\def\e{{\mathbb E}}
\def\p{{\mathbb P}}
\def\P{{\bf P}}
\def\E{{\bf E}}
\def\Q{{\bf Q}}
\def\z{{\mathbb Z}}
\def\N{{\mathbb N}}
\def\T{{\mathbb T}}
\def\G{{\mathbb G}}
\def\L{{\mathbb L}}

\def\deg{\chi}

\def\ee{\mathrm{e}}
\def\d{\, \mathrm{d}}
\def\S{\mathscr{S}}
\def\bs{{\tt bs}}



\vglue50pt

\centerline{\Large\bf Asymptotics for the survival probability}

\bigskip

\centerline{\Large\bf in a killed branching random walk}

\bigskip
\bigskip

\centerline{by}

\medskip

\centerline{Nina Gantert, $\;$Yueyun Hu $\;$and$\;$ Zhan Shi}

\medskip

\centerline{\it Universit\"at M\"unster, Universit\'e Paris XIII \& Universit\'e Paris VI}

\bigskip

\centerline{\tt This version: February 15, 2010}

\bigskip
\bigskip

{\leftskip=2truecm
\rightskip=2truecm
\baselineskip=15pt
\small

\noindent{\slshape\bfseries Summary.} Consider a discrete-time
one-dimensional supercritical branching random walk. We study the
probability that there exists an infinite ray  in the branching
random walk that always lies above the line of slope
$\gamma-\varepsilon$, where $\gamma$ denotes the asymptotic speed
of the right-most position in the branching random walk. Under
mild general assumptions upon the distribution of the branching
random walk, we prove that when $\varepsilon\to 0$, this
probability decays like $\exp\{ - {\beta + o(1)\over
\varepsilon^{1/2}}\}$, where $\beta$ is a positive constant
depending on the distribution of the branching random walk. In the
special case of i.i.d.\ Bernoulli$(p)$ random variables (with
$0<p<{1\over 2}$) assigned on a rooted binary tree, this answers
an open question of Robin Pemantle, see \cite{pemantle}.

\bigskip

\noindent{\slshape\bfseries Keywords.} Branching random walk, survival probability, maximal displacement.

\bigskip

\noindent{\slshape\bfseries 2000 Mathematics Subject Classification.} 60J80.

} 

\bigskip
\bigskip

\section{Introduction}
   \label{s:intro}

We consider a one-dimensional branching random walk in discrete time. Before introducing the model and the problem, we start with an example, borrowed from Pemantle~\cite{pemantle}, in the study of binary search trees.

\medskip

\begin{example}
 \label{ex:pemantle}
 {\rm
Let $\T_\bs$ be a binary tree (``{\tt bs}" for binary search),
rooted at $e$. Let $(Y(x), \, x\in \T_\bs)$ be a collection,
indexed by the vertices of the tree, of i.i.d.\ Bernoulli random
variables with mean $p\in (0, \, {1\over 2})$. For any vertex
$x\in \T_\bs \backslash \{ e\}$, let $[\![ e, \, x]\!]$ denote the
shortest path connecting $e$ with $x$, and let $]\!] e, \, x]\!]
:= [\![ e, \, x]\!] \backslash \{ e\}$. We define
$$
U_\bs (x) := \sum_{v\in \, ]\!] e, \, x]\!]} Y(v), \qquad x\in \T_\bs \backslash \{ e\},
$$

\noindent and $U_\bs(e) := 0$. Then $(U_\bs(x), \, x\in \T_\bs )$ is a binary branching Bernoulli random walk. It is known (Kingman~\cite{kingman}, Hammersley~\cite{hammersley}, Biggins~\cite{biggins}) that
$$
\lim_{n\to \infty} \, {1\over n} \max_{|x|=n} U_\bs(x) = \gamma_\bs, \qquad \mbox{a.s.,}
$$

\noindent where the constant $\gamma_\bs=\gamma_\bs(p)\in (0, \, 1)$ is the unique solution of
\begin{equation}
    \gamma_\bs\log {\gamma_\bs\over p} +
    (1-\gamma_\bs) \log {1-\gamma_\bs\over 1-p}
    - \log 2=0.
    \label{c(p)}
\end{equation}

For any $\varepsilon>0$, let $\varrho_\bs(\varepsilon,p)$ denote
the probability that there exists an infinite ray\footnote{By an
infinite ray, we mean that each $x_j$ is the parent of $x_{j+1}$.}
$\{ e=: x_0, \, x_1, \, x_2, \, \ldots \}$ such that $U_\bs(x_j)
\ge (\gamma_\bs-\varepsilon)j$ for all $j\ge 1$. It is conjectured
by Pemantle~\cite{pemantle} that there exists a constant
$\beta_\bs(p)$ such that\footnote{Throughout the paper, by
$a(\varepsilon) \sim b(\varepsilon)$, $\varepsilon\to 0$, we mean
$\lim_{\varepsilon\to 0} {a(\varepsilon) \over b(\varepsilon)}
=1$.}
\begin{equation}
    \log \varrho_\bs(\varepsilon,p)
    \; \sim \; - {\beta_\bs(p) \over
    \varepsilon^{1/2}},
    \qquad \varepsilon\to 0.
    \label{pemantle}
\end{equation}

We prove the conjecture, and give the value of $\beta_\bs(p)$. Let $\psi_\bs(t) := \log [2(p \ee^t + 1-p)]$, $t>0$. Let $t^* = t^*(p) >0$ be the unique solution of $\psi_\bs(t^*) = t^* \psi_\bs '(t^*)$. [One can then check that the solution of equation (\ref{c(p)}) is $\gamma_\bs = {\psi_\bs (t^*)\over t^*}$.] Our main result, Theorem \ref{t:main} below, implies that conjecture (\ref{pemantle}) holds, with
$$
\beta_\bs (p) := {\pi \over 2^{1/2}} [t^* \psi_\bs''(t^*)]^{1/2} .
$$

A particular value of $\beta_\bs$ is as follows: if $0<p_0<{1\over 2}$ is such that $16p_0(1-p_0)=1$ (i.e., if $\gamma_\bs(p_0)= {1\over 2}$), then
$$
\beta_\bs (p_0) = {\pi\over 4} \Big( {\gamma_\bs'(p_0) \over 1-2p_0}\Big)^{1/2} \log {1\over 4p_0} ,
$$

\noindent where $\gamma_\bs'(p_0)$ denotes the derivative  of
$p\mapsto \gamma_\bs(p)$ at $p_0$. This is, informally, in
agreement with the following theorem of Aldous (\cite{aldous},
Theorem 6): if $p\in (p_0, \, {1\over 2})$ is such that
$\gamma_\bs (p) = {1\over 2} + \varepsilon$, then the probability
that there exists an infinite ray $x$ with $U_\bs (x_i) \ge
{1\over 2} \, i$, $\forall i\ge 1$, is
$$
\exp\Big( - {\pi \log (1/(4p_0)) \over 4(1-2p_0)^{1/2}} \, {1\over (p-p_0)^{1/2}} + O(1) \Big), \qquad \varepsilon\to 0.
\eqno\sqcup\!\!\!\!\sqcap
$$

} 
\end{example}

\bigskip

As a matter of fact, the main result of this paper (Theorem \ref{t:main} below) is valid for more general branching random walks: the tree $\T_\bs$ can be random (Galton--Watson), the random variables assigned on the vertices of the tree are not necessarily Bernoulli, nor necessarily identically distributed, nor necessarily independent if the vertices share a common parent.

Our model is as follows, which is a one-dimensional discrete-time branching random walk. At the beginning, there is a single particle located at position $x=0$. Its children, who form the first generation, are positioned according to a certain point process. Each of the particles in the first generation gives birth to new particles that are positioned (with respect to their birth places) according to the same point process; they form the second generation. The system goes on according to the same mechanism. We assume that for any $n$, each particle at generation $n$ produces new particles independently of each other and of everything up to the $n$-th generation.

We denote by $(U(x), \, |x|=n)$ the positions of the particles in the $n$-th generation, and by $Z_n := \sum_{|x|=n} 1$ the number of particles in the $n$-th generation. Clearly, $(Z_n, \, n\ge 0)$ forms a Galton--Watson process. [In Example \ref{ex:pemantle}, $Z_n= 2^n$, whereas $(U(x), \, |x|=1)$ is a pair of independent Bernoulli$(p)$ random variables.]

We assume that for some $\delta>0$,
\begin{equation}
    \E( Z_1^{1+\delta})<\infty, \qquad
    \E(Z_1)>1;
    \label{delta}
\end{equation}

\noindent in particular, the Galton--Watson process $(Z_n, \, n\ge 0)$ is supercritical. We also assume that there exist $\delta_+>0$ and $\delta_->0$ such that
\begin{equation}
    \E\Big( \sum_{|x|=1} \ee^{\delta_+ U(x)}
    \Big)
    <\infty,
    \qquad
    \E\Big( \sum_{|x|=1} \ee^{-\delta_- U(x)}
    \Big)
    <\infty.    
    \label{delta+-}
\end{equation}

An additional assumption is needed (which in Example \ref{ex:pemantle} corresponds to the condition $p<{1\over 2}$). Let us define the logarithmic generating function for the branching walk:
\begin{equation}
    \psi(t) := \log \E \Big( \sum_{|x|=1}
    \ee^{t U(x)} \Big), \qquad t>0.
    \label{psi}
\end{equation}

\noindent Let $\zeta := \sup \{ t: \psi(t) <\infty\}$. Under
Condition (\ref{delta+-}), we have $0< \zeta \le \infty$, and
$\psi$ is $C^\infty$ on $(0, \, \zeta)$. We assume that there
exists $t^* \in (0, \, \zeta)$ such that
\begin{equation}
    \psi(t^*) = t^* \psi'(t^*) .
    \label{t*}
\end{equation}

\noindent For discussions on this condition, see the examples presented after Theorem \ref{t:main} below.

Recall that (Kingman~\cite{kingman}, Hammersley~\cite{hammersley}, Biggins~\cite{biggins}) conditioned on the survival of the system,
\begin{equation}
    \lim_{n\to \infty} \,
    {1\over n} \max_{|x|=n} U(x) = \gamma,
    \qquad \mbox{a.s.,}
    \label{LGN}
\end{equation}

\noindent where $\gamma := {\psi(t^*)\over t^*}$ is a constant, with $t^*$ and $\psi(\cdot)$ defined in (\ref{t*}) and (\ref{psi}), respectively.

For $\varepsilon>0$, let $\varrho_U(\varepsilon)$ denote the
probability that there exists an infinite ray $\{ e=: x_0, \, x_1,
\, x_2, \, \ldots \}$ such that $U(x_j) \ge (\gamma-\varepsilon)j$
for all $j\ge 1$. Our main result is as follows.

\medskip

\begin{theorem}
 \label{t:main}
 Assume $(\ref{delta})$ and $(\ref{delta+-})$.
 If $(\ref{t*})$ holds, then
 \begin{equation}
     \log \varrho_U(\varepsilon) \; \sim \;
     - {\pi\over (2\varepsilon)^{1/2}}
     \, [t^* \psi''(t^*)]^{1/2},
     \qquad \varepsilon \to 0,
     \label{main}
 \end{equation}
 where $t^*$ and $\psi$ are as in
 $(\ref{t*})$ and $(\ref{psi})$, respectively.
\end{theorem}

\medskip

Since $(U(x), \; |x|=1)$ is not a deterministic set (excluded by the combination of (\ref{t*}) and (\ref{delta})), the function $\psi$ is strictly convex on $(0, \, \zeta)$. In particular, we have $0<\psi''(t^*)<\infty$.

We now present a few simple examples to illustrate the meaning of Assumption (\ref{t*}). For more detailed discussions, see Jaffuel~\cite{jaffuel}.

\bigskip

\noindent {\bf Example \ref{ex:pemantle}} (continuation). In Example \ref{ex:pemantle}, Conditions $(\ref{delta})$ and $(\ref{delta+-})$ are obviously satisfied, whereas $(\ref{t*})$ is equivalent to $p<{1\over 2}$. In this case, (\ref{main}) becomes (\ref{pemantle}). Clearly, if $p>{1\over 2}$, $\varrho_\bs(\varepsilon,p)$ does not go to $0$ because the vertices labeled with $1$ percolate, with positive probability, on the tree. \hfill$\Box$

\medskip

\begin{example}
 \label{ex:bernoulli}
 {\rm
 Consider the example of Bernoulli branching
 random walk, i.e., such that
 $U(x) \in \{0, \, 1\}$ for any $|x|=1$; to
 avoid trivial cases, we assume
 $\E(\sum_{|x|=1} {\bf 1}_{\{ U(x) =1\} })>0$
 and
 $\E(\sum_{|x|=1} {\bf 1}_{\{ U(x) =0\} })>0$.

 Condition $(\ref{delta+-})$ is automatically
 satisfied as long as we assume
 $(\ref{delta})$. Elementary computations
 show that Condition $(\ref{t*})$ is
 equivalent to
 $\E(\sum_{|x|=1} {\bf 1}_{\{ U(x) =1\} })<1$. (In
 particular, if we assign independent
 Bernoulli$(p)$ random variables on the vertices of a
 rooted binary tree, we recover Example \ref{ex:pemantle}).
 Again, if  $\E(\sum_{|x|=1} {\bf 1}_{\{ U(x) =1\} })> 1$, $\varrho_U(\varepsilon)$ does not go to $0$ because the vertices labeled with $1$ percolate, with positive probability, on the tree. \hfill$\Box$

 }

\end{example}

\begin{example}
 \label{ex:borne}
 {\rm
 Assume the distribution of $U$ is bounded from
 above, in the sense that there exists a 
 constant
 $C\in \r$ such that $\sup_{|x|=1}U(x)\le C$.
 Let $s_U :=
 \hbox{\rm ess}\sup \sup_{|x|=1}U(x)
 = \sup\{a\in \r: \, \P\{\sup_{|x|=1}U(x)\ge a\}>0\}
 <\infty$. Under $(\ref{delta})$ and
 $(\ref{delta+-})$, Condition $(\ref{t*})$ is
 satisfied if and only if
 $\E(\sum_{|x|=1} {\bf 1}_{ \{ U(x) = s_U\} })
 <1$.\hfill$\Box$

} 

\end{example}

\begin{example}
 \label{ex:unbound}
 {\rm
 Assume that $(\ref{delta})$ holds true.
 If $\hbox{\rm ess}\sup \sup_{|x|=1}U(x) =
 \infty$, then Condition $(\ref{t*})$ is satisfied. \hfill$\Box$

} 

\end{example}

\medskip

We mention that the question we address here in the discrete case has a continuous counterpart, which has been investigated in the context of the F-KPP equation with cutoff, see \cite{brunet-derrida}, \cite{derrida-simon}, \cite{derrida-simon2}.

\medskip

The rest of the paper is organized as follows. In Section \ref{s:reduction}, we make a linear transformation of our branching random walk so that it will become a boundary case in the sense of Biggins and Kyprianou~\cite{biggins-kyprianou}; the linear transformation is possible due to Assumption $(\ref{t*})$. Section \ref{s:ub} is devoted to the proof of the upper bound in Theorem \ref{t:main}, whereas the proof of the lower bound is in Section \ref{s:lb}.

\section{A linear transformation}
\label{s:reduction}

We define
\begin{equation}
    V(x) := -t^* U(x) + \psi(t^*) \, |x|.
    \label{V}
\end{equation}

\noindent Then
\begin{equation}
    \E\Big( \sum_{|x|=1} \ee^{-V(x)} \Big)=1,
    \qquad
    \E\Big( \sum_{|x|=1} V(x) \ee^{-V(x)}
    \Big)= 0.
    \label{biggins-kyprianou}
\end{equation}

\noindent Since $t^*< \zeta$, there exists $\delta_1>0$ such that 
\begin{equation}
    \E \Big( 
    \sum_{|x|=1} \ee^{-(1+\delta_1)V(x)}
    \Big) <\infty .
    \label{delta1}
\end{equation} 

\noindent On the other hand, by (\ref{delta+-}), there exists $\delta_2>0$ such that
\begin{equation}
    \E \Big( 
    \sum_{|x|=1} \ee^{\delta_2 V(x)}
    \Big) <\infty .
    \label{delta2}
\end{equation} 

The new branching random walk $(V(x))$ satisfies $\lim_{n\to \infty} {1\over n} \inf_{|x|=n} V(x) =0$ a.s.\ conditioned on non-extinction. Let
\begin{equation}
    \varrho(\varepsilon)=
    \varrho(V,\varepsilon) :=
    \P \Big\{ \hbox{\rm $\exists$ infinite ray
    $\{ e=: x_0, \, x_1, \, x_2, \, \ldots \}$:
    $V(x_j) \le \varepsilon j$, $\forall j\ge 1$}
    \Big\}.
    \label{rho}
\end{equation}

\noindent Theorem \ref{t:main} will be a consequence of the following estimate: assuming (\ref{biggins-kyprianou}), then
\begin{equation}
    \log \varrho(\varepsilon) \; \sim \;
    - {\pi \sigma\over (2\varepsilon)^{1/2}},
    \qquad \varepsilon \to 0,
    \label{main2}
\end{equation}

\noindent where $\sigma$ is the constant in (\ref{sigma}) below.

It is (\ref{main2}) we are going to prove: an upper bound is proved in Section \ref{s:ub}, and a lower bound in Section \ref{s:lb}.

We conclude this section with a change-of-probabilities formula, which is the raison d'\^etre of the linear transformation. Let $S_0:=0$, and let $(S_i-S_{i-1}, \; i\ge 1)$ be a sequence of i.i.d.\ random variables such that for any measurable function $f: \, \r\to [0, \, \infty)$,
$$
\E(f(S_1) ) = \E\Big( \sum_{|x|=1} \ee^{-V(x)} f(V(x)) \Big) .
$$

\noindent In particular, $\E(S_1)=0$ (by (\ref{biggins-kyprianou})). In words, $(S_n)$ is a mean-zero random walk. We denote
\begin{equation}
    \sigma^2 := \E(S_1^2)
    = \E \Big(\sum_{|x|=1} V(x)^2 \ee^{-V(x)}
    \Big) = (t^*)^2 \psi''(t^*).
    \label{sigma}
\end{equation}

\noindent Since $\E(Z_1^{1+\delta})<\infty$ (Condition (\ref{delta})) and $\E(\sum_{|x|=1} \ee^{-(1+\delta_1)V(x)}) <\infty$ (see (\ref{delta1})), there exists $\delta_3>0$ such that $\E(\ee^{u S_1}) <\infty$ for all $|u|\le \delta_3$. 

In view of (\ref{biggins-kyprianou}), we have, according to Biggins and Kyprianou~\cite{biggins-kyprianou}, for any $n\ge 1$ and any measurable function $F: \, \r^n \to [0, \, \infty)$,
\begin{equation}
    \E\Big( \sum_{|x|=n} \ee^{-V(x)}
    F(V(x_i), \; 1\le i\le n) \Big)
    = \E [ F(S_i, \; 1\le i\le n) ],
    \label{change-proba}
\end{equation}

\noindent where, for any $x$ with $|x|=n$, $\{ e=: x_0, \, x_1, \, \dots, x_n:=x\}$ is the shortest path connecting $e$ to $x$. 

We now give a bivariate version of (\ref{change-proba}). For any vertex $x$, the number of its children is denoted by $\nu(x)$. Condition (\ref{delta}) guarantees that $\P\{ \nu(x)<\infty, \; \forall x\} =1$. In light of (\ref{biggins-kyprianou}), we have, for any $n\ge 1$ and any measurable function $F: \, \r^n\times \r^n \to [0, \, \infty)$,
\begin{equation}
    \E\Big( \sum_{|x|=n} \ee^{-V(x)}
    F[V(x_i), \; 
    \nu(x_{i-1}), \; 1\le i\le n] \Big)
    = \E \Big( 
    F[S_i, \; \nu_{i-1}, \; 1\le i\le n]
    \Big),
    \label{change-proba2}
\end{equation}

\noindent where $(S_i-S_{i-1}, \; \nu_{i-1})$, for $i\ge 1$, are i.i.d.\ random vectors, whose common distribution is determined by (recalling that $Z_1 := \# \{ y: \; |y|=1\}$)
\begin{equation}
    \E[f(S_1, \; \nu_0)] 
    = 
    \E\Big( \sum_{|x|=1} \ee^{-V(x)} 
    f ( V(x), \; Z_1 ) \Big) ,
    \label{nu}
\end{equation}

\noindent for any measurable function $f: \, \r^2\to [0, \, \infty)$.

The proof of (\ref{change-proba2}), just as the proof of (\ref{change-proba}) in Biggins and Kyprianou~\cite{biggins-kyprianou}, relies on a simple argument by induction on $n$. We feel free to omit it. 

[We mention that (\ref{change-proba2}) is a special case of the so-called spinal decomposition for branching random walks, a powerful tool developed by Lyons, Pemantle and Peres~\cite{lyons-pemantle-peres} and Lyons~\cite{lyons}. The idea of spinal decomposition, which goes back at least to Kahane and Peyri\`ere~\cite{kahane-peyriere}, has been used in the literature by many authors in several different forms.]

We now extend a useful result of Mogulskii~\cite{mogulskii} to {\it arrays} of random variables.

\medskip

\begin{lemma}  
\label{l:mogulskii-array}
 {\bf (A triangular version of 
 Mogulskii~\cite{mogulskii})}
 For each $n\ge 1$, let
 $X_i^{(n)}$, $1\le i \le n$, be i.i.d.\
 real-valued random variables.
 Let $g_1<g_2$ be continuous functions on
 $[0, \, 1]$ with $g_1(0) < 0 < g_2(0)$.
 Let $(a_n)$ be a sequence of positive numbers 
 such that 
 $a_n \to\infty$ and that 
 ${a^2_n\over n} \to 0$.
 Assume that there exist constants 
 $\eta>0$ and $\sigma^2>0$ such that
 \begin{equation}
     \sup_{n\ge 1} 
     \E (|X_1^{(n)}|^{2+\eta}) <\infty,
     \quad  
     \E(X_1^{(n)})=o\Big( {a_n\over n} \Big),
     \quad 
     \mbox{\rm Var}(X_1^{(n)}) \to \sigma^2.
     \label{mogulskii-condition}
 \end{equation}
 Consider the measurable event
 $$
 E_n := \Big\{ g_1\Big({i\over n}\Big) \le
 {S_i^{(n)}\over a_n} \le 
 g_2\Big({i\over n}\Big), \;
 \hbox{ \it for } 1\le i\le n \Big\} ,
 $$
 where $S_i^{(n)} := X_1^{(n)} + \cdots +
 X_i^{(n)}$. We have
 \begin{equation}
     \lim_{n\to \infty} {a_n^2\over n} \log
     \P \{E_n \} =
     - {\pi^2 \sigma^2\over 2} \,
     \int_0^1 {\d t\over [g_2(t)- g_1(t)]^2}.
     \label{mogulskii}
 \end{equation}
 Moreover, for any $b>0$,
 \begin{equation}
     \lim_{n\to \infty} {a_n^2\over n} \log
     \P\Big\{ E_n, \; {S_n^{(n)}\over a_n} \ge
     g_2(1)-b\Big\} =
     - {\pi^2 \sigma^2\over 2} \,
     \int_0^1 {\d t\over [g_2(t)- g_1(t)]^2} .
     \label{mogulskii2}
 \end{equation}

\end{lemma}

\medskip

If the distribution of $X_1^{(n)}$ does not depend on $n$, Lemma \ref{l:mogulskii-array} is Mogulskii~\cite{mogulskii}'s result. In this case, condition (\ref{mogulskii-condition}) is satisfied as long as $X_1^{(n)}$ is centered, having a finite $(2+\eta)$-moment (for some $\eta>0$), and such that it is not identically zero.\footnote{In this case, we even can allow $\eta=0$; see \cite{mogulskii}.}

The proof of Lemma \ref{l:mogulskii-array} is in the same spirit (but with some additional technical difficulties) as in the original work of Mogulskii~\cite{mogulskii}, and is included as an appendix at the end of the paper. We mention that as in \cite{mogulskii}, it is possible to have a version of Lemma \ref{l:mogulskii-array} when $X^{(n)}_1$ belongs to the domain of attraction of a stable non-Gaussian law, except that the constant ${\pi^2\over 2}$ in (\ref{mogulskii})--(\ref{mogulskii2}) will be implicit.

\section{Proof of Theorem \ref{t:main}: the upper bound}
\label{s:ub}

In this section, we prove the upper bound in (\ref{main2}):
\begin{equation}
    \limsup_{\varepsilon\to 0}
    \varepsilon^{1/2} \log \varrho(\varepsilon)
    \le - {\pi \sigma\over 2^{1/2}} \, ,
    \label{ub}
\end{equation}

\noindent where $\varrho(\varepsilon)$ is defined in (\ref{rho}), and $\sigma$ is the constant in (\ref{sigma}).

The main idea in this section is borrowed from Kesten~\cite{kesten}. We start with the trivial inequality that for any $n\ge 1$ (an appropriate value for $n= n(\varepsilon)$ will be chosen later on),
$$
\varrho(\varepsilon) \le \P\Big\{ \exists \, x: \; |x|=n, \; V(x_i) \le \varepsilon i, \; \forall i\le n \Big\}.
$$

\noindent Let $(b_i, \, i\ge 0)$ be a sequence of non-negative real numbers whose value (depending on $n$) will be given later on. For any $x$, let $H(x) := \inf\{ i: \, 1\le i\le |x|, \, V(x_i) \le \varepsilon i - b_i\}$, with $\inf \emptyset := \infty$. Then $\P \{ H(x) =\infty\} + \P\{ H(x) \le |x|\} =1$. Therefore,
$$
\varrho (\varepsilon) \le \varrho_1 (\varepsilon)  + \varrho_2 (\varepsilon),
$$

\noindent where
\begin{eqnarray*}
    \varrho_1 (\varepsilon)
    = \varrho_1 (\varepsilon, n)
 &:=& \P\Big\{ \exists \, |x|=n: \;
    H(x) = \infty , \; V(x_i) \le \varepsilon i,
    \; \forall i\le n \Big\} ,
    \\
    \varrho_2 (\varepsilon)
    = \varrho_2 (\varepsilon, n)
 &:=& \P\Big\{ \exists \, |x|=n: \; H(x) \le n,
    \; V(x_i) \le \varepsilon i, \;
    \forall i\le n \Big\}.
\end{eqnarray*}

\noindent We now estimate $\varrho_1 (\varepsilon)$ and $\varrho_2 (\varepsilon)$ separately. 

By definition,
\begin{eqnarray*}
    \varrho_1 (\varepsilon)
 &=& \P\Big\{ \exists \, |x|=n: \;
    \varepsilon i - b_i<V(x_i) \le \varepsilon i,
    \; \forall i\le n \Big\}
    \\
 &=& \P\Big\{ \sum_{|x|=n}
    {\bf 1}_{\{ \varepsilon i - b_i< V(x_i) \le
    \varepsilon i,
    \; \forall i\le n\} } \ge 1 \Big\}
    \\
 &\le& \E \Big( \sum_{|x|=n}
    {\bf 1}_{\{ \varepsilon i - b_i < V(x_i) \le
    \varepsilon i, \; \forall i\le n\} } \Big),
\end{eqnarray*}

\noindent the last inequality being a consequence of Chebyshev's
inequality. Applying the change-of-probabilities formula
(\ref{change-proba}) to $F(z) := \ee^{z_n} {\bf 1}_{\{ \varepsilon
i - b_i < z_i \le \varepsilon i, \; \forall i\le n\} }$ for $z:=
(z_1, \ldots , z_n) \in \r^n$, this yields, in the notation of
(\ref{change-proba}),
\begin{equation}
    \varrho_1 (\varepsilon) \le
    \E \Big( \ee^{S_n}
    {\bf 1}_{\{ \varepsilon i - b_i < S_i \le
    \varepsilon i, \; \forall i\le n\} } \Big)
    \le
    \ee^{\varepsilon n}\,
    \P\Big\{ \varepsilon i - b_i< S_i \le
    \varepsilon i, \; \forall i\le n \Big\}.
    \label{rho1}
\end{equation}

To estimate $\varrho_2(\varepsilon)$, we observe that
\begin{eqnarray*}
    \varrho_2(\varepsilon)
 &\le & \sum_{j=1}^n \P\Big\{ \exists \, |x|=n: \;
    H(x) =j, \; V(x_i) \le \varepsilon i, \;
    \forall i\le n \Big\}
    \\
 &\le& \sum_{j=1}^n \P\Big\{ \exists \, |x|=n: \;
    H(x) =j, \; V(x_i) \le \varepsilon i, \;
    \forall i\le j \Big\} .
\end{eqnarray*}

\noindent Since $\{ \exists \, |x|=n: \; H(x) =j, \; V(x_i) \le \varepsilon i, \; \forall i\le j \} \subset \{ \exists \, |y|=j: \; H(y) =j, \; V(y_i) \le \varepsilon i, \; \forall i\le j \}$, this yields
$$
\varrho_2(\varepsilon) \le \sum_{j=1}^n \P\Big\{ \exists \, |y|=j: \; \varepsilon i -b_i < V(y_i) \le \varepsilon i, \; \forall i<j, \; V(y_j) \le \varepsilon j -b_j \Big\} .
$$

\noindent We can now use the same argument as for $\varrho_1(\varepsilon)$, namely, Chebyshev's inequality and then the change-of-probability formula (\ref{biggins-kyprianou}), to see that
\begin{eqnarray*}
    \varrho_2(\varepsilon)
 &\le& \sum_{j=1}^n \E\Big( \sum_{|y|=j}
    {\bf 1}_{ \{ \varepsilon i -b_i < V(y_i) \le
    \varepsilon i, \; \forall i<j, \;
    V(y_j) \le \varepsilon j -b_j\} } \Big)
    \\
 &=& \sum_{j=1}^n \E\Big( \ee^{S_j}
    {\bf 1}_{ \{ \varepsilon i -b_i < S_i \le
    \varepsilon i, \; \forall i<j, \;
    S_j \le \varepsilon j -b_j\} } \Big)
    \\
 &\le& \sum_{j=1}^n \ee^{\varepsilon j -b_j}
    \P\Big\{ \varepsilon i -b_i < S_i \le
    \varepsilon i, \; \forall i<j \Big\}.
\end{eqnarray*}

\noindent Together with (\ref{rho1}), and recalling that $\varrho(\varepsilon) \le \varrho_1(\varepsilon) + \varrho_2(\varepsilon)$, this yields
\begin{eqnarray*}
    \varrho (\varepsilon)
 &\le& \ee^{\varepsilon n}\,
    \P\Big\{ \varepsilon i - b_i< S_i \le
    \varepsilon i, \; \forall i\le n \Big\} +
    \sum_{j=1}^n \ee^{\varepsilon j -b_j}
    \P\Big\{ \varepsilon i -b_i < S_i \le
    \varepsilon i, \; \forall i<j \Big\}
    \\
 &=&\ee^{\varepsilon n} I(n) +
    \sum_{j=0}^{n-1}
    \ee^{\varepsilon (j+1) -b_{j+1}} I(j),
\end{eqnarray*}

\noindent where $I(0):=1$ and
$$
I(j) := \P\Big\{ \varepsilon i -b_i < S_i \le \varepsilon i, \; \forall i\le j \Big\} , \qquad 1\le j\le n.
$$

The idea is now to apply Mogulskii's estimate (\ref{mogulskii}) to $I(j)$ for suitably chosen $(b_i)$. Unfortunately, since $\varepsilon$ depends on $n$, we are not allowed to apply (\ref{mogulskii}) simultaneously to all $I(j)$, $0\le j\le n$. So let us first work a little bit more, and then apply (\ref{mogulskii}) to only a few of the $I(j)$.

We assume that $(b_i)$ is non-increasing. Fix an integer $N\ge 2$, and take $n:= kN$ for $k\ge 1$. Then
\begin{eqnarray}
    \varrho(\varepsilon)
 &\le& \ee^{\varepsilon kN} I(kN) +
    \sum_{j=0}^{k-1}
    \ee^{\varepsilon (j+1) -b_{j+1}} I(j) +
    \sum_{\ell =1}^{N-1}
    \sum_{j= {\ell k}}^{{(\ell+1)k}-1}
    \ee^{\varepsilon (j+1) -b_{j+1}} I(j)
    \nonumber
    \\
 &\le& \ee^{\varepsilon kN} I(kN) +
    k\exp ( \varepsilon k -b_k) +
    k\sum_{\ell =1}^{N-1}
    \exp \Big( {\varepsilon (\ell +1)k} -
    b_{(\ell +1)k} \Big) I(\ell k).
    \label{rho<}
\end{eqnarray}

We choose $b_i = b_i(n) := b (n-i)^{1/3}= b (kN-i)^{1/3}$, $0\le i\le n$, and $\varepsilon := {\theta \over n^{2/3}}= {\theta \over (Nk)^{2/3}}$, where $b>0$ and $\theta>0$ are constants. By definition, for $1\le \ell \le N$,
$$
I(\ell k) =  \P\Big\{ \theta \Big({\ell\over N}\Big)^{2/3} {i\over \ell k}
- b \Big({N\over \ell}-{i\over \ell k}\Big)^{1/3} < {S_i \over (\ell
k)^{1/3}} \le \theta \Big({\ell\over N}\Big)^{2/3} {i\over \ell k}, \;
\forall i \le \ell k\Big\}.
$$

\noindent Applying (\ref{mogulskii}) to $g_1(t) := \theta ({\ell\over N})^{2/3} t - b({N\over \ell} - t)^{1/3}$ and $g_2(t) := \theta ({\ell\over N})^{2/3} t$, we see that, for $1\le \ell \le N$,
$$
\limsup_{k\to \infty} {1\over (\ell k)^{1/3}} \log I(\ell k) \le - {\pi^2 \sigma^2 \over 2b^2} \int_0^1 {\d t\over ({N\over \ell} - t)^{2/3}} = - {3\pi^2 \sigma^2 \over 2b^2} \, {N^{1/3} - (N-\ell)^{1/3} \over \ell^{1/3}},
$$

\noindent where $\sigma$ is the constant in (\ref{sigma}). Going back to (\ref{rho<}), we obtain:
$$
\limsup_{k\to \infty} {\theta^{1/2} \over (Nk)^{1/3}} \log \varrho\Big( {\theta \over (Nk)^{2/3}} \Big) \le \theta^{1/2} \alpha_{N,b},
$$

\noindent where the constant $\alpha_{N,b}= \alpha_{N,b}(\theta)$ is defined by
\begin{eqnarray*}
    \alpha_{N,b}
 &:=& \max_{1\le \ell \le N-1}
    \Big\{ \theta -{3\pi^2 \sigma^2 \over 2b^2} ,
    \; {\theta\over N} - b(1- {1\over N})^{1/3},
    \\
 && \qquad\qquad {\theta (\ell +1)\over N} -
    b(1- {\ell+1\over N})^{1/3} -
    {3\pi^2 \sigma^2 \over 2b^2} \, {N^{1/3} -
    (N-\ell)^{1/3} \over N^{1/3}} \Big\}.
\end{eqnarray*}

\noindent Since $\varepsilon\mapsto \varrho(\varepsilon)$ is non-increasing, this yields
$$
\limsup_{\varepsilon \to 0} \varepsilon^{1/2} \log  \varrho(\varepsilon) \le \theta^{1/2} \alpha_{N,b}.
$$

We let $N\to \infty$. By definition,
$$
\limsup_{N\to \infty} \alpha_{N,b}
\le \max \Big\{ \theta - {3\pi^2 \sigma^2 \over 2b^2} , \; - b, \; f(\theta,b) \Big\} ,
$$

\noindent where $f(\theta,b) := \sup_{t\in (0, \, 1]} \{ \theta t- b(1- t)^{1/3} - {3\pi^2 \sigma^2 \over 2b^2} [1- (1-t)^{1/3}] \}$.

Elementary computations show that as long as $b<{3\pi^2 \sigma^2
\over 2b^2} \le  b+ 3\theta$, we have $f(\theta,b) = \theta -
{3\pi^2 \sigma^2 \over 2b^2} + {2\over 3 (3\theta)^{1/2}} ({3\pi^2
\sigma^2 \over 2b^2} -b)^{3/2}$.  Thus $\max \{ \theta - {3\pi^2
\sigma^2 \over 2b^2} , \; - b, \; f(\theta,b) \} = \max \{
f(\theta,b) , \; - b\}$, which equals $-b$ if $\theta = {\pi^2
\sigma^2 \over 2b^2} - {b\over 3}$. As a consequence, for any
$b>0$   satisfying $b < {3\pi^2 \sigma^2 \over 2b^2}$,
$$
\limsup_{\varepsilon \to 0} \varepsilon^{1/2} \log
\varrho(\varepsilon) \le - b  \, \sqrt{ {\pi^2 \sigma^2 \over
2b^2} - {b\over 3}} = - \sqrt{ {\pi^2 \sigma^2 \over 2 } -
{b^3\over 3}}.
$$

\noindent Letting $b\to 0$, this   yields (\ref{ub})  and
completes the proof of the upper bound in Theorem
\ref{t:main}.\hfill$\Box$

\section{Proof of Theorem \ref{t:main}: the lower bound}
\label{s:lb}

Before proceeding to the proof of the lower bound in Theorem \ref{t:main}, we recall two inequalities: the first gives a useful lower tail estimate for the number of individuals in a super-critical Galton--Watson process conditioned on survival, whereas the second concerns an elementary property of the conditional distribution of a sum of independent random variables. Let us recall that $Z_n$ is the number of particles in the $n$-th generation.

\medskip

\begin{fact}
 \label{f:mcdiarmid}
 {\bf (McDiarmid~\cite{mcdiarmid})}
 There exists $\vartheta>1$ such that
 \begin{equation}
     \P \{ Z_n \le \vartheta^n
     \, | \, Z_n >0 \}
     \le \vartheta^{-n},
     \qquad \forall n\ge 1.
     \label{mcdiarmid}
 \end{equation}

\end{fact}

\medskip

\begin{fact}
 \label{f:yzpolymer}
 {\bf (\cite{yzpolymer})}
 If $X_1, X_2, \ldots ,X_N$ are
 independent non-negative random variables,
 and if $F: (0, \, \infty) \to\r_+$ is
 non-increasing, then
 $$
 \E \Big[ F\Big( \sum_{i=1}^N X_i\Big) \, \Big| \,
 \sum_{i=1}^N X_i >0 \Big]
 \le \max_{1\le i\le N}
 \E [F(X_i) \, | \, X_i >0] .
 $$

\end{fact}

\medskip

This section is devoted to the proof of the lower bound in (\ref{main2}):
\begin{equation}
    \liminf_{\varepsilon\to 0}
    \varepsilon^{1/2} \log \varrho(\varepsilon)
    \ge - {\pi \sigma\over 2^{1/2}} \, ,
    \label{lb}
\end{equation}

\noindent where $\varrho(\varepsilon)$ and $\sigma$ are as in (\ref{rho}) and (\ref{sigma}), respectively.

The basic idea consists in constructing a new Galton--Watson tree $\G=\G(\varepsilon)$ within the branching random walk, and obtaining a lower bound for $\varrho(\varepsilon)$ in terms of $\G$.

Recall from (\ref{LGN}) that conditioned on survival, ${1\over j} \max_{|z|\le j} V(z)$ converges almost surely, for $j\to \infty$, to a finite constant. [The fact that this limiting constant is finite is a consequence of $\E(\sum_{|x|=1} \ee^{\delta_2 V(x)})<\infty$ in (\ref{delta2}).] Since the system survives with (strictly) positive probability, we can fix a sufficiently large constant $M>0$ such that
\begin{equation}
    \inf_{j\ge 0}
    \P \Big\{ \max_{|x|\le j} V(x) \le M j 
    \Big\}
    \ge {1\over 2}, \qquad
    \kappa:= \inf_{j\ge 0}
    \P \Big\{ Z_j>0, \;
    \max_{|x|\le j} V(x) \le M j \Big\} >0,
    \label{M}
\end{equation}

\noindent where, as before, $Z_j := \# \{x: |x|=j\}$.

Fix a constant $0<\alpha<1$. For any integers $n>L\ge 1$ with $(1-\alpha)\varepsilon L \ge M (n-L)$, we consider the set $G_{n,\varepsilon}= G_{n,\varepsilon}(L)$ defined by
\footnote{We write $z>x$ if $x$ is an ancestor of $z$.}
$$
G_{n,\varepsilon}:= \{|x|= n: \, V(x_i) \le \alpha \varepsilon i,
 \;  \hbox{\rm for} \; 1\le i\le L; \; \max_{z>x_L:\, |z|\le n}
 [V(z)-V(x_L)] \le (1-\alpha) \varepsilon L \} .
$$

\noindent By definition, for any $x\in G_{n,\varepsilon}$, we have $V(x_i) \le \varepsilonÊi$, for $1\le i\le n$.

If $G_{n,\varepsilon} \not= \emptyset$, the elements of $G_{n,\varepsilon}$ form the first generation of the new Galton--Watson tree $\G_{n,\varepsilon}$, and we construct $\G_{n,\varepsilon}$ by iterating the same procedure: for example, the second generation in $\G_{n,\varepsilon}$ consists of $y$ with $|y|=2n$ being a descendant of some $x\in G_{n,\varepsilon}$ such that $V(y_{n+i})-V(x) \le \alpha \varepsilon i$, for $1 \le i\le L$ and $\max_{z>y_{n+L}:\, |z|\le 2n}[ V(z) -V(y_{n+L})] \le (1-\alpha) \varepsilon L$.

Let $q_{n,\varepsilon}$ denote the probability of extinction of the Galton--Watson tree $\G_{n,\varepsilon}$. It is clear that
$$
\varrho(\varepsilon) \ge 1-q_{n,\varepsilon},
$$

\noindent so we only need to find a lower bound for $1-q_{n,\varepsilon}$. In order to do so, we introduce, for $b\in \r$ and $n\ge 1$,
 \begin{equation}
    \varrho(b, \, n) :=
    \P\Big\{ \exists |x|=n: \; 
    V(x_i) \le b i, 
    \hbox{ \rm for } 1\le i\le n\Big\}.
    \label{rho3}
 \end{equation}

\noindent Let us first prove some preliminary results.

\medskip

\begin{lemma}
 \label{l:lb1}
 Let $0<\alpha<1$ and $\varepsilon>0$. Let
 $n>L\ge 1$ be such that
 $(1-\alpha)\varepsilon L \ge M (n-L)$. Then
 $$
 \P\{ G_{n,\varepsilon} \not= \emptyset \} \ge
 {1\over 2} \varrho(\alpha\varepsilon, \, n).
 $$
\end{lemma}

\medskip

\noindent {\it Proof.} By definition,
$$
\P\{ G_{n,\varepsilon} \not= \emptyset \} = \E \Big( {\bf 1}_{ \{ \exists |y|=L: \; V(y_i) \le \alpha \varepsilon i, \; \forall i\le L\} }
 \P \Big\{ \max_{|z|\le n-L} V(z) \le (1-\alpha)\varepsilon L \Big\} \Big).
$$

\noindent Since $(1-\alpha)\varepsilon L \ge M (n-L)$, it follows from (\ref{M}) that
$$
\P\{ G_{n,\varepsilon} \not= \emptyset \} \ge {1\over 2} \,
\P\{ \exists |y|=L: \; V(y_i) \le \alpha \varepsilon i, \; \forall i\le L \},
$$

\noindent and the r.h.s. is at least  ${1\over 2} \varrho(\alpha\varepsilon, \, n)$.\hfill$\Box$

\medskip

\begin{lemma}
 \label{l:lb2}
 Let $0<\alpha<1$ and $\varepsilon>0$. Let
 $n>L\ge 1$ be such that
 $(1-\alpha)\varepsilon L \ge M (n-L)$. We have
 \begin{equation}\label{Gest}
 \P\{ 1\le \# G_{n,\varepsilon} \le \vartheta^{n-L} \} \le
 {1\over \kappa \, \vartheta^{n-L}},
 \end{equation}
 where $\kappa>0$ and $\vartheta>1$ are the
 constants in $(\ref{M})$ and
 $(\ref{mcdiarmid})$, respectively.
\end{lemma}

\medskip

\noindent {\it Proof.} By definition,
$$
\# G_{n,\varepsilon} = \sum_{|x|=L} \eta_x {\bf 1}_{ \{ V(x_i) \le \alpha \varepsilon i, \; \forall i\le L\} } ,
$$

\noindent where
$$
\eta_x := \# \{ y>x: \; |y|=n\} \, {\bf 1}_{ \{ \max_{\{ z>x: \, |z| \le n\} } [V(z)-V(x)] \le (1-\alpha) \varepsilonÊL \} }.
$$

\noindent By Fact \ref{f:yzpolymer}, for any $\ell \ge 1$, with $F(x) = {\bf 1}_{\{ x \leq \ell\} }$,
$$
\P\Big\{ \# G_{n,\varepsilon} \le \ell \, \Big| \, \# G_{n,\varepsilon} >0\Big\} \le \P\Big\{ Z_{n-L} \le \ell \, \Big| \, Z_{n-L}>0, \; \max_{|z| \le n-L} V(z) \le (1-\alpha) \varepsilonÊL \Big\},
$$

\noindent where, as before, $Z_{n-L} := \# \{ |x| = n-L\}$. Since $(1-\alpha) \varepsilonÊL \ge M (n-L)$, it follows from (\ref{M}) that
$\P \{Z_{n-L}>0, \; \max_{|z| \le n-L} V(z) \le (1-\alpha) \varepsilonÊL \} \ge \kappa>0$. Therefore,
$$
\P \{1\le \# G_{n,\varepsilon} \le \ell \} \le {1\over \kappa} \, \P \Big\{ Z_{n-L} \le \ell \, \Big| \, Z_{n-L}>0 \Big\}.
$$

\noindent This implies (\ref{Gest}) by means of Fact \ref{f:mcdiarmid}.\hfill$\Box$

\bigskip

To state the next estimate, we recall that $\nu(x)$ is the number of children of $x$, and that $(S_i-S_{i-1}, \, \nu_{i-1})$, $i\ge 1$, are i.i.d.\ random vectors (with $S_0:=0$) whose common distribution is given by $(\ref{nu})$.

\medskip

\begin{lemma}
 \label{l:lb8}
 Let $n\ge 1$. For any $1\le i\le n$, let 
 $I_{i,n}\subset \r$ be a Borel set. Let 
 $r_n \ge 1$ be an integer. We have
 $$
 \P\Big\{ \exists |x|=n : \;
 V(x_i) \in I_{i,n} \, , \;
 \forall 1\le i\le n\Big \}
 \ge {\E [\ee^{S_n} \, 
 {\bf 1}_{\{ S_i \in I_{i,n}\, , \; 
 \nu_{i-1} \le r_n, \; \forall \, 
 1\le i\le n\} } ]
 \over 1+ (r_n-1)\sum_{j=1}^n h_{j,n}} ,
 $$
 where 
 \begin{equation}
     h_{j,n}:= 
     \sup_{u\in I_{j,n}} 
     \E \Big( \ee^{S_{n-j}}
     {\bf 1}_{ \{ S_\ell \in 
     I_{\ell +j,n}-u, \; 
     \forall 0\le \ell\le n-j\} } \Big) ,
     \label{h}
 \end{equation}
 and 
 $I_{\ell +j,n}-u 
 := \{ v-u: \; v\in I_{\ell +j,n} \}$.

\end{lemma}

\medskip

\noindent {\it Proof.} Let
$$
Y_n := \sum_{|x|=n} 
{\bf 1}_{\{ V(x_i) \in I_{i,n} \, , \;
\nu (x_{i-1}) \le r_n, 
\; \forall 1\le i\le n\} }.
$$

\noindent By definition,
\begin{eqnarray}
    \E(Y_n^2)
 &=& \E \Big(
    \sum_{|x|=n}\sum_{|y|=n}
    {\bf 1}_{ \{ V(x_i)\in I_{i,n} \, , \;
    \nu(x_{i-1}) \le r_n, \;
    V(y_i) \in I_{i,n} \, , \; 
    \nu(y_{i-1}) \le r_n , \;
    \forall 1\le i\le n\} } \Big)
    \nonumber
    \\
 &=& \E(Y_n) + 
    \E \Big( \sum_{j=0}^{n-1} \sum_{|z|=j}
    {\bf 1}_{ \{ V(z_i) \in I_{i,n} \, , \;
    \nu(z_{i-1}) \le r_n, \;
    \forall i\le j\} } D_{j+1,n}(z) \Big) ,
    \label{E(Yn2)}
\end{eqnarray}

\noindent with
\begin{eqnarray*}
    D_{j+1,n}(z)
 &:=& \sum_{(x_{j+1}, \, y_{j+1})}
    \sum_{(x, \, y)}
    {\bf 1}_{ \{ V(x_i) \in I_{i,n} \, ,\;
    \nu(x_{i-1}) \le r_n, \;
    V(y_i) \in I_{i,n} \, , \; 
    \nu(y_{i-1}) \le r_n , \;
    \forall j+1 \le i\le n\} } 
    \\
 &\le& \sum_{(x_{j+1}, \, y_{j+1})}
    \sum_{(x, \, y)}
    {\bf 1}_{ \{ V(x_i) \in I_{i,n} \, ,\;
    \nu(x_{i-1}) \le r_n, \;
    V(y_i) \in I_{i,n} \, , \;
    \forall j+1 \le i\le n\} } ,
\end{eqnarray*} 

\noindent where the double sum $\sum_{(x_{j+1}, \, y_{j+1})}$ is over pairs $(x_{j+1},\, y_{j+1})$ of distinct children of $z$ (thus $|x_{j+1}|=|y_{j+1}|=j+1$), while $\sum_{(x,\, y)}$ is over pairs $(x,\, y)$ with $|x|=|y|=n$ such that\footnote{We write $y\ge x$ if either $y>x$ or $y=x$.} $x\ge x_{j+1}$ and $y\ge y_{j+1}$. 

The $\E[\sum_{j=0}^{n-1} \sum_{|z|=j} {\bf 1}_{ \{ \cdots \} } D_{j+1,n}(z)]$ expression on the right-hand side of (\ref{E(Yn2)}) is bounded by
$$
\E \Big( \sum_{j=0}^{n-1} \sum_{|z|=j} {\bf 1}_{ \{ V(z_i) \in I_{i,n} \, , \; \nu(z_{i-1}) \le r_n, \; \forall i\le j\} } \sum_{(x_{j+1}, \, y_{j+1})} \sum_x {\bf 1}_{ \{ V(x_i) \in I_{i,n} \, , \; \nu(x_{i-1}) \le r_n, \; \forall j+1\le i\le n\} } \, h_{j+1,n} \Big),
$$

\noindent where $h_{j+1,n}:= \sup_{u\in I_{j+1,n}} \E [ \sum_{|y|=n-j-1}{\bf 1}_{ \{ V(y_\ell) \in I_{\ell +j+1,n}-u, \; \forall 0\le \ell\le n-j-1\} }]$, which is in agreement with (\ref{h}), thanks to the change of probability formula (\ref{change-proba}). [The sum $\sum_x$ is, of course, still over $x$ with $|x|=n$ such that $x\ge x_{j+1}$.]

Thanks to the condition $\nu(x_j) \le r_n$ (i.e., $\nu(z) \le r_n$), we see that the sum $\sum_{y_{j+1}}$ in the last display gives at most a factor of $r_n-1$; which yields that the last display is at most $(r_n-1) \E(\sum_{j=0}^{n-1} Y_n h_{j+1,n})$. In other words, we have proved that
$$
\E \Big( \sum_{j=0}^{n-1} \sum_{|z|=j}
    {\bf 1}_{ \{ V(z_i) \in I_{i,n} \, ,\; 
    \nu(z_{i-1}) \le r_n, \;
    \forall i\le j\} } D_{j+1,n}(z) \Big) \le (r_n-1) \sum_{j=0}^{n-1} \E(Y_n) h_{j+1,n} .
$$

\noindent This yields $\E(Y_n^2) \le [1+ (r_n-1)\sum_{j=0}^{n-1} h_{j+1,n}]\E(Y_n)$. Therefore, 
\begin{equation}
    {\E(Y_n^2) \over [\E(Y_n)]^2} 
    \le 
    {1+ (r_n-1)\sum_{j=1}^n h_{j,n} \over  
    \E(Y_n)} 
    = 
    {1+ (r_n-1)\sum_{j=1}^n h_{j,n} \over
    \E ( \ee^{S_n} \, 
    {\bf 1}_{\{ S_i \in I_{i,n}\, , \; 
    \nu_{i-1} \le r_n, \; 
    \forall 1\le i\le n\} } )} ,
    \label{paley-zygmund}
\end{equation}

\noindent the last inequality being a consequence of (\ref{change-proba2}). By the Cauchy--Schwarz inequality, $\P \{ Y_n \ge 1 \} \ge {[\E(Y_n)]^2 \over \E(Y_n^2)}$. Recalling the definition of $Y_n$, we obtain from (\ref{paley-zygmund}) that
\begin{eqnarray}
 &&\P\Big\{ \exists |x|=n : \;
    V(x_i) \in I_{i,n} \, , \;
    \nu (x_{i-1}) \le r_n, 
    \; \forall 1\le i\le n\Big \} 
    \nonumber
    \\
 &\ge& {\E [\ee^{S_n} \, 
    {\bf 1}_{\{ S_i \in I_{i,n}\, , \; 
    \nu_{i-1} \le r_n, \; \forall \, 
    1\le i\le n\} } ]
    \over 1+ (r_n-1)\sum_{j=1}^n h_{j,n}} .
    \label{lb8}
\end{eqnarray}
    
\noindent Lemma \ref{l:lb8} follows immediately from (\ref{lb8}).\hfill$\Box$

\bigskip

The key step in the proof of the lower bound in Theorem \ref{t:main} is the following estimate.

\medskip

\begin{lemma}
 \label{l:lb4}
 For any $\theta>0$,
 $$
 \liminf_{n\to \infty}
 {\log \varrho(\theta n^{-2/3}, \, n)\over
 n^{1/3}} \ge
 -{\pi \sigma \over (2\theta)^{1/2}},
 $$
 where $\sigma>0$ is the constant in
 $(\ref{sigma})$.
\end{lemma}

\medskip

\noindent {\it Proof.} Let $0<\lambda<{\pi \sigma \over (2\theta)^{1/2}}$, and let $I_{i,n} := [{\theta i\over n^{2/3}}- \lambda \, n^{1/3}, \; {\theta i\over n^{2/3}}]$ (for $1\le i\le n$). Since $\varrho(\theta n^{-2/3}, \, n) \ge \P\{ \exists |x|=n: \; V(x_i) \in I_{i,n} \, , \; \forall 1\le i\le n\}$, it follows from Lemma \ref{l:lb8} that for any integer $r_n \ge 1$,
$$
\varrho(\theta n^{-2/3}, \, n) \ge {\E [\ee^{S_n} \, 
 {\bf 1}_{\{ S_i \in I_{i,n}\, , \; 
 \nu_{i-1} \le r_n, \; \forall \, 
 1\le i\le n\} } ]
 \over 1+ (r_n-1)\sum_{j=1}^n h_{j,n}}
=: {\Lambda_n \over 1+ (r_n-1)\sum_{j=1}^n h_{j,n}} ,
$$

\noindent where $h_{j,n}$ is defined in $(\ref{h})$, while $(S_i-S_{i-1}, \, \nu_{i-1})$, $i\ge 1$, are i.i.d.\ random vectors (with $S_0:=0$) whose common distribution is given by $(\ref{nu})$.

For any $\theta_1<\theta$, we have
\begin{eqnarray*}
    \Lambda_n
 &\ge& \ee^{\theta_1 n^{1/3}} \,
    \P\{ S_i \in I_{i,n} \, , \;
    \nu_{i-1} \le r_n , \;
    \forall 1\le i\le n, \;
    S_n \ge \theta_1 n^{1/3}\}
    \\
 &=& \ee^{\theta_1 n^{1/3}} \,
    \P\Big\{ \theta {i\over n}
    -\lambda \le {S_i\over n^{1/3}}
    \le \theta {i\over n}, \;
    \nu_{i-1} \le r_n , \;
    \forall 1\le i\le n, \;
    {S_n \over n^{1/3}} \ge \theta_1\Big\}.
\end{eqnarray*}

For any $n\ge 1$, we consider i.i.d.\ random variables $X_i^{(n)}$, $1\le i\le n$, having the same distribution as $S_1$ conditioned on $\nu_0 \le r_n$. Let $S_0^{(n)} =0$ and $S_i := X_1^{(n)} + \cdots + X_i^{(n)}$ for $1\le i\le n$. Then
$$
\Lambda_n \ge \ee^{\theta_1 n^{1/3}} \, [ \P\{\nu_0 \le r_n\}]^n \, \P\Big\{ \theta {i\over n} -\lambda \le {S_i^{(n)}\over n^{1/3}} \le \theta {i\over n}, \; \forall 1\le i\le n, \; {S_n^{(n)} \over n^{1/3}} \ge \theta_1\Big\} .
$$

We now choose $r_n := \lfloor \ee^{n^{1/4}}\rfloor$. By definition, $\P\{\nu_0 > r_n\} = \E (\sum_{|x|=1} \ee^{-V(x)} {\bf 1}_{ \{ Z_1 > r_n\} })$, where $Z_1 = \sum_{|y|=1} 1$ as before. By Markov's inequality, $\P\{ Z_1 > r_n\} \le {\E(Z_1^{1+\delta})\over r_n^{1+\delta}}$. Since $\E(Z_1^{1+\delta})<\infty$ (Condition (\ref{delta})) and $\E(\sum_{|x|=1} \ee^{-(1+\delta_1) V(x)})<\infty$ (see (\ref{delta1})), an application of H\"older's inequality confirms that $\P\{\nu_0 > r_n\} \le r_n^{-\delta_4}$ for some $\delta_4>0$ and all sufficiently large $n$. In view of our choice of $r_n$, we see that $[\P\{\nu_0 \le r_n\}]^n \to 1$. Therefore, for all sufficiently large $n$,
$$
\Lambda_n \ge {1\over 2} \, \ee^{\theta_1 n^{1/3}} \, \P\Big\{ \theta {i\over n} -\lambda \le {S_i^{(n)}\over n^{1/3}} \le \theta {i\over n}, \; \forall 1\le i\le n, \; {S_n^{(n)} \over n^{1/3}} \ge \theta_1\Big\} .
$$

To deal with the probability expression on the right-hand side, we intend to apply (\ref{mogulskii2}); so we need to check condition (\ref{mogulskii-condition}). Recall that $S_1$ has finite exponential moments in the neighbourhood of 0. Thus, the first condition in (\ref{mogulskii-condition}), namely, $\sup_{n\ge 1} \E (|X_1^{(n)}|^{2+\eta}) <\infty$ for some $\eta>0$, is trivially satisfied. To check the second condition, we see that since $\E(S_1)=0$, we have $\E(X_1^{(n)}) = - {\E[S_1\, {\bf 1}_{ \{ \nu_0> r_n\} }] \over \P\{\nu_0 \le r_n\}}$. Since $\P\{\nu_0 > r_n\} \le r_n^{-\delta_4}$ for some $\delta_4>0$ and all sufficiently large $n$, and since $S_1$ has some finite exponential moments, the second condition in (\ref{mogulskii-condition}), $\E(X_1^{(n)})=o ( {a_n\over n})$, is also satisfied (regardless of the value of the sequence $a_n \to \infty$) in view of the Cauchy--Schwarz inequality. Moreover, $\E(X_1^{(n)}) \to 0$, which yields $\mbox{\rm Var}(X_1^{(n)}) \to \E(S_1^2) - 0 = \sigma^2$: the third and last condition in (\ref{mogulskii-condition}) is verified.

We are therefore entitled to apply (\ref{mogulskii2}): taking $g_1(t) := \theta t-\lambda$ and $g_2(t) := \theta t$, we see that for any $\lambda_1 \in (0, \, \lambda)$ and all sufficiently large $n$,
$$
    \Lambda_n
    \ge
    {1\over 2}\,
    \ee^{\theta_1 n^{1/3}}
    \exp\Big( -
    {\pi^2 \sigma^2 \over 2\lambda_1^2}\, n^{1/3}
    \Big),
$$

\noindent which implies, for all sufficiently large $n$,
\begin{equation}
    \varrho(\theta n^{-2/3}, \, n) \ge
    {{1\over 2} 
    \exp[ (\theta_1 -
    {\pi^2 \sigma^2 \over 2\lambda_1^2})
    n^{1/3}]
    \over 1+ (r_n-1)\sum_{j=1}^n h_{j,n}} .
    \label{rho>}
\end{equation}

To estimate $\sum_{j=1}^n h_{j,n}$, we observe that 
\begin{eqnarray*}
    h_{j,n}
 &=& \sup_{u\in I_{j,n}}
    \E \Big( \ee^{S_{n-j}}
    {\bf 1}_{ \{ S_i \in
    [{\theta(i+j)\over n^{2/3}}-\lambda n^{1/3}
    -u, \;
    {\theta(i+j)\over n^{2/3}}-u], \;
    \forall 0\le i\le n-j\} } \Big)
    \\
 &=& \sup_{v\in [0, \; \lambda n^{1/3}]}
    \E \Big( \ee^{S_{n-j}}
    {\bf 1}_{ \{ S_i \in
    [{\theta i\over n^{2/3}}-\lambda n^{1/3}
    +v, \;
    {\theta i\over n^{2/3}}+v], \;
    \forall 0\le i\le n-j\} } \Big)
    \\
 &\le& \ee^{\theta(n-j)n^{-2/3}+\lambda n^{1/3}}
    \sup_{v\in [0, \; \lambda n^{1/3}]}
    \P\Big\{ {\theta i\over n^{2/3}}-\lambda n^{1/3}
    +v
    \le S_i \le {\theta i\over n^{2/3}} +v, \;
    \forall 0\le i\le n-j\Big\} .
\end{eqnarray*}

We now use the same trick as in the proof of the upper bound in Theorem \ref{t:main} by sending $n$ to infinity along a subsequence. Fix an integer $N\ge 1$. Let $n:=Nk$, with $k\ge 1$. For any $j\in [(\ell-1)k+1, \, \ell k]\cap \z$ (with $1\le \ell \le N$), we have
$$
h_{j,n}\le \ee^{\theta(N-\ell +1)kn^{-2/3} +\lambda n^{1/3}} \sup_{v\in [0, \; \lambda n^{1/3}]} \P\Big\{ v-\lambda n^{1/3} \le S_i - {\theta i\over n^{2/3}} \le v, \; \forall i\le (N-\ell)k\Big\}.
$$

\noindent Unfortunately, the interval $[0, \, \lambda n^{1/3}]$ in $\sup_{v\in [0, \; \lambda n^{1/3}]} \P\{\cdots\}$ is very large, so we split it into smaller ones of type $[{(m-1)\lambda n^{1/3}\over N}, \; {m\lambda n^{1/3}\over N}]$ (for $1\le m\le N$), to see that the $\sup_{v\in [0, \; \lambda n^{1/3}]} \P\{\cdots\}$ expression is
\begin{eqnarray*}
 &\le& \max_{1\le m\le N}
    \P\Big\{ {(m-1)\lambda n^{1/3}\over N}
    -\lambda n^{1/3}
    \le S_i - {\theta i\over n^{2/3}} \le
    {m\lambda n^{1/3}\over N}, \;
    \forall i\le (N-\ell)k\Big\}
    \\
 &=&\max_{1\le m\le N}
    \P\Big\{ - {(N-m+1)\lambda \over N^{2/3}} \le {S_i\over k^{1/3}}- {\theta \over N^{2/3}} {i\over k} \le
    {m\lambda \over N^{2/3}}, \;
    \forall i\le (N-\ell)k\Big\}.
\end{eqnarray*}

\noindent We are now entitled to apply (\ref{mogulskii}) to $n:=(N-\ell)k$, $g_1(t) := {\theta \over (N-\ell)^{1/3} N^{2/3}} t - {(N-m+1)\lambda \over (N-\ell)^{1/3} N^{2/3}}$ and $g_2(t) := {\theta \over (N-\ell)^{1/3} N^{2/3}} t + {m\lambda \over (N-\ell)^{1/3} N^{2/3}}$, to see that for any $1\le \ell \le N$ and uniformly in $j\in [(\ell-1)k+1, \, \ell k]\cap \z$ (and in $j=0$, which formally corresponds to $\ell=0$),
$$
\limsup_{k\to \infty} {1\over N^{1/3}k^{1/3}} \log h_{j,Nk} \le {\theta (N-\ell +1) \over N} + \lambda - {\pi^2\sigma^2 \over 2}{(N-\ell) N\over (N+1)^2\lambda^2},
$$

\noindent which is bounded by ${\theta (N+1) \over N} + \lambda - {\pi^2\sigma^2 \over 2}{N^2\over (N+1)^2\lambda^2}$ (recalling that $\theta > {\pi^2 \sigma^2 \over 2\lambda^2}$). As a consequence,
$$
\limsup_{k\to \infty} {1\over N^{1/3}k^{1/3}} \log \sum_{j=0}^n h_{j,Nk} \le {\theta (N+1) \over N} + \lambda - {\pi^2\sigma^2 \over 2}{N^2\over (N+1)^2\lambda^2} =: c(\theta,N,\lambda).
$$

\noindent Going back to (\ref{rho>}), we get
$$
\liminf_{k\to \infty} {\log \varrho(\theta N^{-2/3} k^{-2/3}, \, Nk)\over N^{1/3}k^{1/3}} \ge \theta_1 - {\pi^2 \sigma^2\over 2 \lambda_1^2} - c(\theta,N,\lambda) .
$$

\noindent By the monotonicity of $n\mapsto \varrho(\theta n^{-2/3}, \, n)$, we obtain:
$$
\liminf_{n\to \infty} {\log \varrho(\theta n^{-2/3}, \, n)\over n^{1/3}} \ge \theta_1 - {\pi^2 \sigma^2\over 2 \lambda_1^2} - c(\theta,N,\lambda).
$$

\noindent Sending $N\to \infty$, $\theta_1\to \theta$, $\lambda \to {\pi \sigma \over (2\theta)^{1/2}}$ and $\lambda_1\to {\pi \sigma \over (2\theta)^{1/2}}$ (in this order) completes the proof of Lemma \ref{l:lb4}.\hfill$\Box$

\bigskip

We now have all the ingredients for the proof of the lower bound in Theorem \ref{t:main}.

\bigskip

\noindent {\it Proof of Theorem \ref{t:main}: the lower bound.} Fix constants $0<\alpha<1$ and $b>\max\{ {M\over 1-\alpha}, \, {(3\pi \sigma)^2 \over \alpha (\log \vartheta)^2} \}$. Let $n>1$. Let
$$
\varepsilon = \varepsilon(n) := {b\over n^{2/3}}, \qquad L=L(n) := n- \lfloor n^{1/3}\rfloor.
$$

\noindent Then $(1-\alpha)\varepsilon L \ge M (n-L)$ for all sufficiently large $n$, say\footnote{Without further mention, the value of $n_0$ can change from line to line when other conditions are to be satisfied.} $n\ge n_0$.

Consider the moment generating function of the reproduction distribution in the Galton--Watson tree $\G_{n,\varepsilon}$:
$$
f(s) := \E(s^{\# G_{n,\varepsilon}}), \qquad s\in [0, \, 1].
$$

\noindent It is well-known that $q_{n,\varepsilon}$, the extinction probability of $\G_{n,\varepsilon}$, satisfies $q_{n,\varepsilon} = f(q_{n,\varepsilon})$. Therefore, for any $0<r <\min\{ q_{n,\varepsilon}\, , \, {1\over 16}\}$,
$$
q_{n,\varepsilon} = f(0) + \int_0^{q_{n,\varepsilon}} f'(s) \d s  = f(0)+ \int_0^{q_{n,\varepsilon} - r} f'(s) \d s + \int_{q_{n,\varepsilon} -r}^{q_{n,\varepsilon}} f'(s) \d s .
$$

\noindent Since $s\mapsto f'(s)$ is non-decreasing on $[0, \, 1]$, we have $\int_0^{q_{n,\varepsilon} - r} f'(s) \d s \le f'(1-r)$. On the other hand, since $f'(s) \le f'(q_{n,\varepsilon}) \le 1$ for $s\in [0, \, q_{n,\varepsilon}]$, we have $\int_{q_{n,\varepsilon} -r}^{q_{n,\varepsilon}} f'(s) \d s \le r$. Therefore,
$$
q_{n,\varepsilon} \le f(0) + f'(1-r) + r.
$$

\noindent Of course, $f(0) = \P\{ G_{n,\varepsilon} = \emptyset \}$, whereas $f'(1-r) = \E[ (\# G_{n,\varepsilon})  (1-r)^{\# G_{n,\varepsilon} -1}]$, which is bounded by ${1\over 1-r} \E[ (\# G_{n,\varepsilon} ) \ee^{-r \# G_{n,\varepsilon}}]$ (using the elementary inequality $1-u\le \ee^{-u}$ for $u\ge 0$). This leads to (recalling that $r<{1\over 16}< {1\over 2}$):
$$
1- q_{n,\varepsilon} \ge \P\{ G_{n,\varepsilon} \not= \emptyset \} - 2\E[ (\# G_{n,\varepsilon}) \ee^{-r \# G_{n,\varepsilon}}] - r.
$$

\noindent Since $u\mapsto u\ee^{-ru}$ is decreasing on $[{1\over r}, \, \infty)$, we see that $\E[ (\# G_{n,\varepsilon}) \ee^{-r \# G_{n,\varepsilon}}]$ is bounded by $\E[ (\# G_{n,\varepsilon}) \ee^{-r \# G_{n,\varepsilon}}\, {\bf 1}_{\{\# G_{n,\varepsilon} \le r^{-2}\} }] + r^{-2} \ee^{-1/r} \le r^{-2} \P\{1\le \# G_{n,\varepsilon} \le r^{-2}\} + r^{-2} \ee^{-1/r}$. Accordingly,
\begin{eqnarray*}
    1- q_{n,\varepsilon}
 &\ge & \P \{ G_{n,\varepsilon} \not= \emptyset \} -
    {2\over r^2}\, \P\{1\le \# G_{n,\varepsilon} \le r^{-2} \} -
    {2\ee^{-1/r}\over r^2} - r
    \\
 &\ge & {1\over 2} \,
    \varrho(\alpha\varepsilon, \, n) -
    {2\over r^2}\, \P \{1\le \# G_{n,\varepsilon} \le r^{-2} \}
    - 2r,
\end{eqnarray*}

\noindent the last inequality following from Lemma \ref{l:lb1} and the fact that $\sup_{\{0<r\le {1\over 16}\} } {1\over r^3} \ee^{-1/r}  < {1\over 2}$.

We choose $r:= {1\over 16} \, \varrho(\alpha\varepsilon, \, n)$. [Since $\varrho(\varepsilon) \ge 1-q_{n,\varepsilon}$, whereas $\lim_{\varepsilon\to 0} \varrho(\varepsilon)=0$ (proved in Section \ref{s:ub}), we have $q_{n,\varepsilon}\to 1$ for $n\to \infty$, and thus the requirement $0<r <\min\{ q_{n,\varepsilon}\, , \, {1\over 16}\}$ is satisfied for all sufficiently large $n$.]

By Lemma \ref{l:lb4}, $r^{-2} \le \vartheta^{n-L}$ for all $n\ge n_0$ (because ${2\pi \sigma \over (\alpha b)^{1/2}} < \log \vartheta$ by our choice of $b$). Therefore, an application of Lemma \ref{l:lb2} tells us that for $n\ge n_0$, $\P\{1\le \# G_{n,\varepsilon} \le r^{-2}\} \le {1\over \kappa \, \vartheta^{n-L}}$, which, by Lemma \ref{l:lb2} again, is bounded by $r^3$ (because ${3\pi \sigma \over (\alpha b)^{1/2}} < \log \vartheta$). Consequently, for all $n\ge n_0$,
$$
1- q_{n,\varepsilon} \ge {1\over 2} \, \varrho(\alpha\varepsilon, \, n) - 2r - 2r = {1\over 4} \, \varrho(\alpha\varepsilon, \, n).
$$

\noindent Recall that $\varrho(\varepsilon) \ge 1-q_{n,\varepsilon}$. Therefore,
$$
\liminf_{n\to \infty} {1\over n^{1/3}} \log \varrho \Big({b\over n^{2/3}}\Big) \ge -{\pi \sigma \over (2\alpha b)^{1/2}}.
$$

\noindent Since $\varepsilon\mapsto \varrho(\varepsilon)$ is non-increasing, we obtain:
$$
\liminf_{\varepsilon\to 0} \varepsilon^{1/2} \log \varrho(\varepsilon) \ge -{\pi \sigma \over (2\alpha)^{1/2}}.
$$

\noindent Sending $\alpha\to 1$ yields (\ref{lb}), and thus proves the lower bound in Theorem \ref{t:main}.\hfill$\Box$

\section{Appendix. Proof of Lemma \ref{l:mogulskii-array}}
\label{s:appendix}

We write $S^{(n)}_j := \sum_{i=1}^j X^{(n)}_i$ (for $1\le j\le n$) and $S^{(n)}_0 :=0$. We need to prove the lower bound in (\ref{mogulskii2}), and the upper bound in (\ref{mogulskii}).

\bigskip

\noindent {\bf Lower bound in (\ref{mogulskii2}).} We want to prove that for any $b>0$,
$$
     \liminf_{n\to \infty} {a_n^2\over n} \log
     \P\Big\{ E_n, \; {S_n^{(n)}\over a_n} \ge
     g_2(1)-b\Big\} \ge
     - {\pi^2 \sigma^2\over 2} \,
     \int_0^1 {\d t\over [g_2(t)- g_1(t)]^2} .
$$

Let $g: \, [0, \, 1]\to \r$ be a continuous function such that $g_1(t)< g(t) < g_2(t)$ for all $t\in [0, \, 1]$. It suffices to prove the lower bound in (\ref{mogulskii2}) when $b>0$ is sufficiently small; so we assume, without loss of generality, that $g(1) \ge g_2(1)-b$. 

Let $\delta>0$ be such that 
\begin{equation}
    g(t) - g_1(t)>3\delta, \quad
    g_2(t)-g(t) >9\delta, \qquad
    \forall t\in [0, \, 1].
    \label{g}
\end{equation}

\noindent Let $A$ be a sufficiently large integer such that    \begin{equation}
    \sup_{0\le s \le t \le 1: \; t-s \le 
    {2\over A} } ( | g_1(t) - g_1(s) |+ 
    | g (t) -  g (s) |+| g_2(t) - g_2(s) |) 
    \le \delta .
    \label{A}
\end{equation}

\noindent Let $r_n := \lfloor A a_n^2\rfloor$, $N=N(n) := \lfloor {n\over r_n}\rfloor$. Let $m_N := n$ and $m_k := k r_n$ for $0\le k\le N-1$. 

Since $g(1) \ge g_2(1)-b$, we have, by definition,
\begin{eqnarray*}
    \P\Big\{ E_n, \; {S_n^{(n)}\over a_n} \ge
    g_2(1)-b\Big\} 
 &\ge& \P\Big(\bigcap_{k=1}^N 
    \Big\{ g_1\big({i\over n}\big) \le
    {S^{(n)}_i\over a_n} \le 
    g_2\big({i\over n}\big), \;
    \forall i\in (m_{k-1}, \, m_k] \cap \z, \;
    \\
 &&\qquad\qquad g \big({m_k\over n}\big) \le
     {S^{(n)}_{m_k}\over a_n} \le
     g \big({m_k\over n}\big) 
     + 6\delta\Big\}\Big) .
\end{eqnarray*} 

\noindent Applying the Markov property successively at times $m_{N-1}$, $m_{N-2}$, $\cdots$, $m_1$, we obtain, by writing $y_k:= g \big({ m_k\over n}\big)$ for $1\le k \le N$,
$$
    \P\Big\{ E_n, \; {S_n^{(n)}\over a_n} \ge
    g_2(1)-b\Big\} 
  \ge p_{1,n}(0) \times \prod_{k=2}^N 
    \inf_{y\in [y_{k-1},\, y_{k-1} + 6\delta]} 
    p_{k,n}(y),
$$

\noindent where\footnote{For notational simplification, we write $\forall i\le \Delta_k$ instead of $\forall i\in (0, \, \Delta_k] \cap \z$.} for $1\le k\le N$ and $y\in \r$,
\begin{eqnarray*} 
    p_{k,n}(y)
 &:=& \P\Big\{ \alpha_{i,k,n} 
    \le {S^{(n)}_i\over a_n} +y \le 
    \beta_{i,k,n}, \;
    \forall i\le \Delta_k; 
    \; y_k  \le 
    {S^{(n)}_{\Delta_k}\over a_n} +y
    \le y_k+ 6\delta \Big\} ,
    \\
    \alpha_{i,k,n}
 &:=& g_1\big({i+m_{k-1}\over n}\big), 
 ÊÊÊ\qquad \beta_{i,k,n} :=
    g_2\big({i+m_{k-1}\over n}\big),
    \qquad
    \Delta_k := m_k - m_{k-1}. 
\end{eqnarray*}

Uniform continuity of $g$ guarantees that when $n$ is sufficiently large, $|y_k-y_{k-1}| \le \delta$ (for all $1\le k\le N$, with $y_0:=0$). In the rest of the proof, we will always assume that $n$ is sufficiently large, say $n\ge n_0$, with $n_0$ depending on $A$ and $\delta$.

We need to bound $p_{1,n}(0) \times \prod_{k=2}^N \inf_{y\in [y_{k-1},\, y_{k-1} + 6\delta]} p_{k,n}(y)$ from below. Let us first get rid of the infimum $\inf_{y\in [y_{k-1},\, y_{k-1} + 6\delta]}$, which is the minimum between $\inf_{y\in [y_{k-1},\, y_{k-1} + 3\delta]}$ and $\inf_{y\in [y_{k-1}+3\delta,\, y_{k-1} + 6\delta]}$: 
$$
\inf_{y\in [y_{k-1},\, y_{k-1} + 6\delta]} p_{k,n}(y) \ge \min\{ p_{k,n}^{(1)}, \, p_{k,n}^{(2)} \}, \qquad n\ge n_0, \; 2\le k\le N,
$$

\noindent where, for $1\le k\le N$,
\begin{eqnarray*}
    p_{k,n}^{(1)}
 &:=& \P\Big\{ 
    \alpha_{i,k,n} - y_{k-1}
    \le {S^{(n)}_i\over a_n} \le 
    \beta_{i,k,n}- y_{k-1} -3\delta , \;
    \forall i\le \Delta_k; 
    \; \delta \le 
    {S^{(n)}_{\Delta_k}\over a_n} 
    \le 2 \delta \Big\}, 
    \\
    p_{k,n}^{(2)}
 &:=& \P\Big\{
    \alpha_{i,k,n} - y_{k-1} -3 \delta
    \le {S^{(n)}_i\over a_n} \le 
    \beta_{i,k,n}- y_{k-1} -6 \delta , \;
    \forall i\le \Delta_k; 
    \; -2\delta \le 
    {S^{(n)}_{\Delta_k}\over a_n} 
    \le -\delta \Big\} .
\end{eqnarray*} 

\noindent And, of course, $p_{1,n}(0) \ge p_{1,n}^{(1)}\ge \min\{ p_{1,n}^{(1)}, \, p_{1,n}^{(2)} \}$. We arrive at the following estimate:
$$
\P\Big\{ E_n, \; {S_n^{(n)}\over a_n} \ge g_2(1)-b\Big\} 
\ge 
\prod_{k=1}^N \min\{ p_{k,n}^{(1)}, \, p_{k,n}^{(2)} \}
= \min\{ p_{N,n}^{(1)}, \, p_{N,n}^{(2)} \}
\prod_{k=1}^{N-1} \min\{ p_{k,n}^{(1)}, \, p_{k,n}^{(2)} \}.
$$

First, we bound $p_{k,n}^{(1)}$ and $p_{k,n}^{(2)}$ from below, for $1\le k\le N-1$ (in which case $\Delta_k = r_n$).
We split the indices $k\in (0, \, N-1] \cap \z$ into $A$ blocs, by means of $(0, \, N-1] \cap \z = \cup_{\ell =1}^A J_\ell$, where $J_\ell = J_\ell(n) := ({(\ell-1)(N-1)\over A}, \, {\ell (N-1)\over A}] \cap \z$. For indices $k$ lying in a same bloc $J_\ell$, we use a common lower bound for $\min\{ p_{k,n}^{(1)}, \, p_{k,n}^{(2)} \}$ as follows: assuming $k\in J_\ell$, we have, by (\ref{A}), $\alpha_{i,k,n} \le g_1({\ell \over A}) + \delta$, $\beta_{i,k,n} \ge g_2({\ell \over A}) - \delta$ (for $n\ge n_0$ and $i\le r_n$), and $|y_{k-1}-g({\ell \over A})| \le \delta$, which leads to: $p_{k,n}^{(1)} \ge q_{\ell,n}^{(1)}$, and $p_{k,n}^{(2)} \ge q_{\ell,n}^{(2)}$ (for $n\ge n_0$, $1\le \ell \le A$ and $k\in J_\ell$), where (recalling that $\Delta_k =r_n$ for $1\le k\le N-1$)
\begin{eqnarray*}
    q_{\ell,n}^{(1)}
 &:=& \P\Big\{ 
    g_1({\ell \over A}) - g({\ell \over A})
     + 2\delta
    \le {S^{(n)}_i\over a_n} \le 
    g_2({\ell \over A}) - g({\ell \over A})
    -5\delta , \;
    \forall i\le r_n; 
    \; \delta \le 
    {S^{(n)}_{r_n}\over a_n} 
    \le 2 \delta \Big\} , 
    \\
    q_{\ell,n}^{(2)} 
 &:=& \P\Big\{ 
    g_1({\ell \over A}) - g({\ell \over A})
     -\delta
    \le {S^{(n)}_i\over a_n} \le 
    g_2({\ell \over A}) - g({\ell \over A})
    -8\delta , \;
    \forall i\le r_n; 
    \; -2\delta \le 
    {S^{(n)}_{r_n}\over a_n} 
    \le -\delta \Big\} .
\end{eqnarray*} 

\noindent Therefore,
$$
    \P\Big\{ E_n, \; {S_n^{(n)}\over a_n} \ge g_2(1)-b\Big\} 
    \ge 
    \min\{ p_{N,n}^{(1)}, \, p_{N,n}^{(2)} \} \, 
    \prod_{\ell=1}^A 
    \Big( \min\{ q_{\ell,n}^{(1)} , \, q_{\ell,n}^{(2)}  \}
    \Big)^{\# J_\ell}.
$$

\noindent Since $\# J_\ell \le {N\over A} \le {n\over r_n A} \le {n\over (Aa_n^2-1)A}$, this yields  
\begin{equation}
    \P\Big\{ E_n, \; {S_n^{(n)}\over a_n} \ge g_2(1)-b\Big\} 
    \ge 
    \min\{ p_{N,n}^{(1)}, \, p_{N,n}^{(2)} \} \, 
    \prod_{\ell=1}^A 
    \Big( \min\{ q_{\ell,n}^{(1)} , \, q_{\ell,n}^{(2)}  \}
    \Big)^{n/[(Aa_n^2-1)A]}.
    \label{lowpen4}
\end{equation}

It is well-known that the linear interpolation function $t (\in [0, 1]) \to { 1 \over a_n}\{ S^{(n)}_{\lfloor r_nt\rfloor} + (r_n t - \lfloor r_nt\rfloor) X^{(n)}_{\lfloor r_nt\rfloor +1}\}$ converges in law to $ (\sigma \sqrt {A}\, W_t, 0\le t \le 1)$, where $W$ denotes a standard one-dimensional Brownian motion.\footnote{In fact, finite-dimensional convergence is easily obtained by verifying Lindeberg's condition in the central limit theorem, whereas tightness is checked using a standard argument, see for example Billingsley~\cite{billingsley}.} So, if we write
$$
f(a,b,c,d) := \P\Big\{ a \le W_t\le b , \; \forall \, t\in [0, \, 1]; \; c \le W_1 \le d\Big\} ,
$$

\noindent for $a<0<b$ and $a\le c<d\le b$, then for any $1\le \ell \le A$, 
\begin{eqnarray*}
    \lim_{n\to \infty} q_{\ell,n}^{(1)}
 &=& f\Big( 
    {g_1({\ell \over A}) - g({\ell \over A}) 
    + 2\delta \over \sigma A^{1/2}}, \; 
    {g_2({\ell \over A}) - g({\ell \over A})
    -5\delta \over \sigma A^{1/2}}, \; 
    {\delta \over \sigma A^{1/2}}, \; 
    {2 \delta\over \sigma A^{1/2}} \Big) ,
    \\
    \lim_{n\to \infty} q_{\ell,n}^{(2)}
 &=& f\Big( 
    {g_1({\ell \over A}) - g({\ell \over A}) 
    -\delta \over \sigma A^{1/2}}, \; 
    {g_2({\ell \over A}) - g({\ell \over A})
    -8\delta \over \sigma A^{1/2}}, \; 
    -{2\delta \over \sigma A^{1/2}}, \; 
    -{\delta\over \sigma A^{1/2}} \Big) .
\end{eqnarray*}

\noindent [Thanks to (\ref{g}), the limits are (strictly) positive.] The function $f$ is explicitly known (see for example, It\^o and McKean~\cite{ito-mckean}, p.~31):
\begin{equation}
    f(a,b,c,d) = 
    \int_c^d {2\over b-a} \sum_{n=1}^\infty 
    \exp\Big( - {n^2 \pi^2 \over 2(b-a)^2}\Big) 
    \sin\Big( {n\pi |a|\over b-a}\Big) 
    \sin\Big( {n\pi (z-a) \over b-a}\Big) \d z,
    \label{ito-mckean}
\end{equation}

\noindent from which it is easily seen that for all $A$ sufficiently large, say $A\ge A_0$ ($A_0$ depending on $\delta$), uniformly in $1\le \ell \le A$,
\begin{eqnarray*}
    \lim_{n\to \infty} q_{\ell,n}^{(1)}
 &\ge& \exp\Big( - {\sigma^2 \pi^2\over 2} 
    {(1+\delta)A \over 
    [g_2({\ell \over A}) - g_1({\ell \over A}) - 
    7\delta]^2}\Big) ,
    \\
    \lim_{n\to \infty} q_{\ell,n}^{(2)}
 &\ge& \exp\Big( - {\sigma^2 \pi^2\over 2} 
    {(1+\delta)A \over 
    [g_2({\ell \over A}) - g_1({\ell \over A}) - 
    7\delta]^2}\Big).
\end{eqnarray*}

\noindent Similarly, we have a lower bound for $\min\{ p_{N,n}^{(1)}, \, p_{N,n}^{(2)}\}$, the only difference being that $\Delta_N$ is not exactly $r_n$ but lies somewhere between $r_n$ and $2r_n$. This time, we only need a rough estimate: there exists a constant $C>0$ such that
$$
\liminf_{n\to \infty} \min\{ p_{N,n}^{(1)}, \, p_{N,n}^{(2)}\} \ge C.
$$

\noindent In view of (\ref{lowpen4}), we get that, for all $A$ sufficiently large (how large depending on $\delta$),
\begin{eqnarray*}  
    \liminf_{n\to\infty} 
    {a_n^2\over n} 
    \log \P\Big\{ E_n, \; 
    {S_n^{(n)}\over a_n} \ge g_2(1)-b\Big\} 
 &\ge& - {\sigma^2 \pi^2\over 2} {1\over A}
    \sum_{\ell=1}^A {1+\delta \over
    [g_2({\ell \over A}) - g_1({\ell \over A}) 
    -7\delta]^2}   
    \\
 &\ge& - {\sigma^2 \pi^2\over 2} (1+2\delta) 
    \int_0^1 {\d t\over 
    [g_2(t) - g_1(t) - 7\delta]^2} .
\end{eqnarray*}

\noindent Letting $A\to \infty$ and $\delta\to 0$ (in this order), we obtain the desired lower bound in (\ref{mogulskii2}).\hfill$\Box$

\bigskip


\noindent {\bf Upper bound in (\ref{mogulskii}).} The upper bound in (\ref{mogulskii}) is needed in this paper only in the form of the original result of Mogulskii~\cite{mogulskii} (i.e., for sequences, instead of arrays, of random variables). We include its proof for the sake of completeness. It is similar to, and easier than, the proof of the lower bound in (\ref{mogulskii2}).

Let $g$ be as before. Let $\delta>0$ and $A>0$ satisfy again (\ref{g}) and (\ref{A}), respectively. Let again $r_n := \lfloor A a_n^2\rfloor$, $N=N(n) := \lfloor {n\over r_n}\rfloor$. Let $m_k := kr_n$ for $0\le k\le N-1$, but we are not interested in $m_N$ any more. Write again $\alpha_{i,k,n} := g_1({i+m_{k-1}\over n})$ and $\beta_{i,k,n} := g_2({i+m_{k-1}\over n})$.

By the Markov property,
$$
\P(E_n) \le \prod_{k=2}^{N-1} \sup_{y\in [g_1({m_{k-1}\over n}), \, g_2({m_{k-1}\over n})]} \widetilde{p}_{k,n}(y),
$$

\noindent where 
$$
\widetilde{p}_{k,n}(y) := 
\P\Big\{ \alpha_{i,k,n} \le {S^{(n)}_i\over a_n} +y \le \beta_{i,k,n}, \; \forall 0 < i\le r_n \Big\} .
$$

\noindent Since $g_1$ and $g_2$ are bounded, we know that $g_1({m_{k-1}\over n})$ and $g_2({m_{k-1}\over n})$ lie in a compact interval, say $[- K \delta, \, K \delta]$ ($K$ being an integer depending on $\delta$). Therefore 
$$ 
\sup_{y\in [g_1({m_{k-1}\over n}), \, g_2({m_{k-1}\over n})]} \widetilde{p}_{k,n}(y) 
\le 
\max_{j\in [-K, \, K-1]\cap \z}
\sup_{y\in [j \delta,\, (j+1) \delta]} \widetilde{p}_{k,n}(y) .
$$

As in the proof of the lower bound in (\ref{mogulskii2}), we cut the interval $(1,N-1]\cap \z$ into $A$ blocs, by means of $(1,N-1]\cap \z = \cup_{\ell =1}^A J_\ell$, where $J_\ell = J_\ell(n) := ( {(\ell-1)(N-2)\over A}+1, \, {\ell(N-2)\over A}+1] \cap \z$. For $k\in J_\ell$, we have, by (\ref{A}), $\alpha_{i,k,n} \ge g_1({\ell \over A}) - \delta$ and $\beta_{i,k,n} \le g_2({\ell \over A}) + \delta$, which leads to: $\sup_{y\in [j \delta,\, (j+1) \delta]}\widetilde{p}_{k,n}(y) \le \widetilde{q}_{\ell,n}(j)$, where
$$
\widetilde{q}_{\ell,n} (j)
:=\P\Big\{ 
    g_1({\ell \over A}) -(j+2)\delta
    \le {S^{(n)}_i\over a_n} \le 
    g_2({\ell \over A}) - (j-1)\delta , \;
    \forall i\le r_n \Big\} . 
$$

\noindent Therefore, 
$$
\P(E_n) \le \prod_{\ell=1}^A \, [\max_{j\in [-K, \, K) \cap \z} \widetilde{q}_{\ell,n}(j)]^{\# J_\ell}.
$$

\noindent We have $\# J_\ell \ge {N-2\over A} -1\ge {n\over A^2 a_n^2} - {3\over A}-1$. On the other hand, for each pair $(\ell, \, j)$, $\widetilde{q}_{\ell,n} (j)$ converges (as $n\to \infty$) to $\P\{ g_1({\ell \over A}) - (j+2)\delta \le \sigmaÊA^{1/2} W_t\le g_1({\ell \over A}) - (j-1)\delta, \; \forall t\in [0, \, 1] \}$, which, in view of (\ref{ito-mckean}), is bounded by $\exp\{ - {\pi^2 \sigma^2 \over 2} {(1-\delta) A \over [g_2(\ell/A) - g_1(\ell/A) - 3\delta]^2}\}$ for all sufficiently large $A$ and uniformly in $(\ell, \, j)$. Accordingly,
\begin{eqnarray*}  
    \limsup_{n\to\infty} {a_n^2\over n} \log \P(E_n) 
 &\le&  - {\pi^2 \sigma^2\over 2} {1\over A} 
    \sum_{\ell=1}^A 
    {1-\delta \over 
    [g_2({\ell \over A}) - g_1({\ell \over A}) - 3\delta]^2}
    \\
 &\le & - {\pi^2 \sigma^2\over 2} (1-2\delta) 
    \int_0^1 {\d t\over [g_2(t) - g_1(t) - 3\delta]^2}, 
\end{eqnarray*} 

\noindent for all sufficiently large $A$. Since $\delta$ can be as close to $0$ as possible, this yields the upper bound in (\ref{mogulskii}).\hfill$\Box$

\bigskip
\bigskip

\noindent {\bf\Large Acknowledgements}

\bigskip

\noindent We are grateful to Jean B\'erard for pointing out two errors in the first draft of the manuscript.

\bigskip
\bigskip


{\footnotesize

\baselineskip=12pt

\noindent
\begin{tabular}{lll}
& Nina Gantert
    & \hskip40pt Yueyun Hu \\
& Fachbereich Mathematik und Informatik
    & \hskip40pt D\'epartement de Math\'ematiques \\
& Universit\"at M\"unster
    & \hskip40pt Universit\'e Paris XIII \\
& Einsteinstrasse 62
    & \hskip40pt 99 avenue J-B Cl\'ement \\
& D-48149 M\"unster
    & \hskip40pt F-93430 Villetaneuse \\
& Germany
    & \hskip40pt France \\
& {\tt gantert@math.uni-muenster.de}
    & \hskip40pt
    {\tt yueyun@math.univ-paris13.fr}
\end{tabular}

\bigskip
\bigskip

Zhan Shi\par Laboratoire de Probabilit\'es UMR 7599\par
Universit\'e Paris VI\par 4 place Jussieu\par F-75252 Paris Cedex
05\par France\par {\tt zhan.shi@upmc.fr}

}


\begin{thebibliography}{99}
\baselineskip=14pt


\bibitem{aldous}
    Aldous, D.~J.\ (1998).
    A Metropolis-type optimization algorithm on
    the infinite tree.
    {\it Algorithmica} {\bf 22}, 388--412.

\bibitem{biggins}
    Biggins, J.~D.\ (1976).
    The first- and last-birth problems for a
    multitype age-dependent branching process.
    {\it Adv.\ Appl.\ Probab.} {\bf 8}, 446--459.

\bibitem{biggins-kyprianou}
    Biggins, J.~D.\ and Kyprianou, A.~E.\ (2005).
    Fixed points of the smoothing transform: the
    boundary case.
    {\it Electron.\ J.\ Probab.} {\bf 10}, Paper
    no.~17, 609--631.

\bibitem{billingsley}
    Billingsley, P.\ (1968).
    {\it Convergence of Probability Measures.}
    Wiley, New York.

\bibitem{brunet-derrida}
    Brunet, \'E.\ and Derrida, B.\ (1997).
    Shift in the velocity of a front due to a cutoff.
    {\it Phys.\ Rev.\ E} {\bf 56}, 
    no.~3, 2597--2604.

\bibitem{derrida-simon}
    Derrida, B.\ and Simon, D.\ (2007).
    The survival probability of a branching random walk in 
    presence of an absorbing wall.
    {\it Europhys.\ Lett.\ EPL} {\bf 78}, 
    no.~6, Paper no.~60006.

\bibitem{derrida-simon2}
    Derrida, B.\ and Simon, D.\ (2008).
    Quasi-stationary regime of a branching random walk in 
    presence of an absorbing wall.
    {\it J.\ Stat.\ Phys} {\bf 131}, 
    no.~2, 203--233.

\bibitem{hammersley}
    Hammersley, J.~M.\ (1974).
    Postulates for subadditive processes.
    {\it Ann.\ Probab.} {\bf 2}, 652--680.

\bibitem{yzpolymer}
    Hu, Y.\ and Shi, Z.\ (2009).
    Minimal position and critical martingale
    convergence in branching random walks,
    and directed polymers on disordered
    trees.
    {\it Ann.\ Probab.} {\bf 37}, 742--789.

\bibitem{ito-mckean}
    It\^o, K.\ and McKean, H.~P.\ (1965).
    {\it Diffusion Processes and their Sample Paths.}
    Springer, Berlin.

\bibitem{jaffuel}
    Jaffuel, B.\ (2009+).
    The critical barrier for the survival of
    the branching random walk with absorption.
    {\tt ArXiv math.PR/0911.2227}

\bibitem{kahane-peyriere}
    Kahane, J.-P.\ and Peyri\`ere, J.\ (1976).
    Sur certaines martingales de Mandelbrot.
    {\it Adv.\ Math.} {\bf 22}, 131--145.

\bibitem{kesten}
    Kesten, H.\ (1978).
    Branching Brownian motion with absorption.
    {\it Stoch.\ Proc.\ Appl.} {\bf 7},
    9--47.

\bibitem{kingman}
    Kingman, J.~F.~C.\ (1975).
    The first birth problem for an age-dependent
    branching process.
    {\it Ann.\ Probab.} {\bf 3}, 790--801.

\bibitem{lyons}
    Lyons, R.\ (1997).
    A simple path to Biggins' martingale
    convergence for branching random walk.
    In: {\it Classical and Modern Branching
    Processes} (Eds.: K.~B.~Athreya and
    P.~Jagers).
    {\it IMA Volumes in Mathematics and its
    Applications} {\bf 84}, 217--221. Springer,
    New York.

\bibitem{lyons-pemantle-peres}
    Lyons, R., Pemantle, R.\ and Peres, Y.\ (1995).
    Conceptual proofs of $L\log L$ criteria for mean
    behavior of branching processes.
    {\it Ann.\ Probab.} {\bf 23}, 1125--1138.

\bibitem{mcdiarmid}
    McDiarmid, C.\ (1995).
    Minimal positions in a branching random walk.
    {\it Ann.\ Appl.\ Probab.} {\bf 5}, 128--139.

\bibitem{mogulskii}
    Mogulskii, A.~A.\ (1974).
    Small deviations in the space of
    trajectories.
    {\it Theory Probab.\ Appl.} {\bf 19},
    726--736.

\bibitem{pemantle}
    Pemantle, R.\ (2009).
    Search cost for a nearly optimal path in a
    binary tree.
    {\it Ann.\ Appl.\ Probab.} {\bf 19},
    1273--1291.

\end{thebibliography}
\end{document}